\documentclass{article}
\usepackage{amsmath,amstext,amsthm,amsfonts}
\usepackage[dvips]{graphicx}

\theoremstyle{plain}
\newtheorem{theorem}{Theorem}[section]
\newtheorem{proposition}[theorem]{Proposition}
\newtheorem{lemma}[theorem]{Lemma}
\newtheorem{corollary}[theorem]{Corollary}
\newtheorem{maintheorem}{Theorem}

\theoremstyle{definition}
\newtheorem{remark}[theorem]{Remark}

\newcommand{\field}[1]{\mathbb{#1}}
\newcommand{\real}{\field{R}}

\renewcommand{\natural}{\field{N}}

\newcommand{\diam}{\operatorname{diam}}

\newcommand{\vep}{\varepsilon}

\newcommand{\SB}{{\cal B}}

\newcommand{\SE}{{\cal E}}

\newcommand{\SL}{{\cal L}}
\newcommand{\SM}{{\cal M}}

\newcommand{\SP}{{\cal P}}
\newcommand{\SQ}{{\cal Q}}

\begin{document}

\title{Existence and uniqueness of maximizing measures \\
for robust classes of local diffeomorphisms}
\author{Krerley Oliveira and Marcelo Viana
\footnote{This work was partially supported by Pronex, CNPq, Fapeal,
and Faperj, Brazil.}}
\date{ }

\maketitle

\begin{abstract}
We prove existence of maximal entropy measures for an open set of non-uniformly expanding local diffeomorphisms
on a compact Riemannian manifold. In this context the topological entropy coincides with the logarithm of the
degree, and these maximizing measures are eigenmeasures of the transfer operator. When the map is topologically
mixing, the maximizing measure is unique and positive on every open set.
\end{abstract}

\section{Introduction}\label{s.introcution}

In its most basic form, the variational principle states that the
topological entropy of a continuous transformation on a compact
space coincides with the supremum of the entropies of the
probability measures invariant under the transformation. We call
\emph{maximizing measure} any invariant probability for which the
supremum is attained. Existence and uniqueness of such measures
has been investigated by many authors, in a wide variety of
situations. However, the global picture is still very much
incomplete.

In this paper we contribute a simple sufficient condition for
existence and uniqueness, applicable to a large class of
transformations. Some examples we have in mind are the
non-uniformly expanding local diffeomorphisms of Alves, Bonatti,
Viana~\cite{ABV00}, which exhibit only positive Lyapunov exponents
at ``most'' points. But our hypothesis, formulated in
\eqref{eq.condicao} below, is a condition of the type that
Buzzi~\cite{Bu99,Bu00} introduced and called entropy-expansivity:
we only ask that the derivative do not expand $k$-dimensional
volume too much, for all $k$ less than the dimension of the
ambient manifold. We show that this implies existence and, if the
transformation is topologically mixing, uniqueness of the
maximizing measure.

\section{Statement of main result}\label{s.statement}

Let $f:M^d \to M^d$ be a $C^1$ local diffeomorphism on a compact
$d$-dimensional Riemannian manifold. Let $p\ge 1$ be the degree of
$f$, that is, the number $\# f^{-1}(x)$ of preimages of any point
$x\in M$. Define
$$
C_k(f) = \max_{x\in M}\|\Lambda^{k}Df(x)\|,
$$
where $\Lambda^k$ represents the $k$th exterior product. We assume
that $f$ satisfies
\begin{equation}\label{eq.condicao}
\max_{1 \leq k \leq d-1} \log C_k(f) < \log p.
\end{equation}
We say that $f:M\to M$ is \emph{topologically mixing} if given any
open set $U$ there exists $N \in \natural$ such that $f^N(U)=M$.
We are going to prove the following

\begin{maintheorem}\label{teoA}
Assume $f$ satisfies \eqref{eq.condicao}. Then $h_{top}(f)=\log
p$, and any maximal eigenmeasure $\mu$ of the transfer operator
$\SL$ is a maximizing measure. In particular, there exists some
maximizing measure for $f$. If $f$ is topologically mixing then
the maximizing measure is unique and positive on open sets.
\end{maintheorem}

The Ruelle-Perron-Frobenius \emph{transfer operator} of $f:M\to M$
is the bounded linear operator $\SL: C(M)\to C(M)$ defined on the
space $C(M)$ of continuous functions $g:M\to\real$ by
$$
\SL g(x) = \sum_{y : f(y)=x} g(y).
$$
Observe that this is a positive operator. Its dual
$\SL^*:\SM(M)\to\SM(M)$ acts on the space of Borel measures of
$M$, by
$$
\int g \, d\SL^*\nu = \int \SL g \, d\nu,
$$
preserving the cone of positive measures, and the subset of
probability measures.

It is easy to see that the spectra of $\SL$ and $\SL^*$ are
contained in the closed disk of radius $p$. We call \emph{maximal
eigenmeasure} any probability measure $\mu$ that satisfies
$$
\SL^*\mu = p \, \mu.
$$
It is well-known that maximal eigenmeasures do exist. A quick
proof goes as follows. Define $G:\SM_1 \to \SM_1$ on the space of
probabilities $\SM_1$ on $M$ by
$$
G(\nu) = \frac{1}{p} \SL^*\nu.
$$
Then $G$ is continuous relative to the weak$^*$ topology on
$\SM_1$. Since $\SM_1$ is a convex compact space, we may use the
Tychonoff-Schauder theorem to conclude that there exists some
probability $\mu$ such $G(\mu)=\mu$. In other words, $\mu$ is a
maximal eigenmeasure. Observe also that $\mu$ is invariant for
$f$. In fact, for every continuous function $g$ we have that $\SL
(g\circ f) (x) = p g(x)$ and
$$
\int (g \circ f) d\mu
 = \frac{1}{p} \int (g\circ f) \, d\SL^*\mu
 = \frac{1}{p} \int \SL(g\circ f) \, d\mu
 = \int g \, d\mu.
$$

The paper is organized as follows. In Section~\ref{s.measures} we
prove that, under our assumptions, any measure with large entropy
has only positive Lyapunov exponents. In
Section~\ref{s.generating} we prove that measures with positive
Lyapunov exponents admits generating partitions with small
diameter. This conclusion uses the notion of hyperbolic times,
that we recall in Section~\ref{s.hyperbolic}. On its turn, it is
used in Section~\ref{s.rokhlin} to show that the entropy of such
measures is given by a simple formula involving the Jacobian.
Using this formula, we prove in Section~\ref{s.existence} that the
topological entropy is $\log p$ and is attained by any maximal
eigenmeasure. Finally, in Section~\ref{s.uniqueness} we prove that
the maximal measure is unique if the transformation is
topologically mixing.

{\bf Acknowledgements.} We are thankful to V\'\i tor Ara\'ujo for
a conversation that helped clarify the arguments in the last section.

\section{Measures with large entropy}\label{s.measures}

By Oseledets~\cite{Os68}, if $\mu$ is an $f$-invariant probability
measure then for $\mu$-almost every point $x\in M$ there is
$k=k(x)\ge 1$, a filtration
$$
T_x M
 = F^1_x \supset \cdots \supset F_x^k \supset F^{k+1}(x) = \{0\},
$$
and numbers
$\hat\lambda_1(x)>\hat\lambda_2(x)>\dots>\hat\lambda_k(x)$ such
$Df(x) F_x^i = F_{f(x)}^i$ and
$$
\hat\lambda_i = \lim_{n\to\infty} \frac{1}{n} \log \|Df^n(x)v\|
$$
for every $v \in F_x^i\setminus F_x^{i+1}$ and $i=1, \ldots, k$.
The numbers $\hat\lambda_i(x)$ are called \emph{Lyapunov
exponents} of $f$ at the point $x$. The \emph{multiplicity} of
$\lambda_i(x)$ is $\dim F_x^i - \dim F_x^{i+1}$. We also write the
Lyapunov exponents as
$$
\lambda_1(x) \geq \lambda_2(x) \geq \dots \geq \lambda_d(x),
$$
where each number is repeated according to the corresponding
multiplicity. Then the \emph{integrated Lyapunov exponents} are
the averages
$$
\lambda_i(\mu)=\int\lambda_i(x)\,d\mu(x),
 \quad\text{for } i=1, \ldots, d.
$$

Given a vector space $V$ and a number $k\ge 1$, the $k$th exterior
power of $V$ is the vector space of all alternate $k$-linear forms
defined on the dual of $V$. We always take $V$ to be
finite-dimensional, and then the exterior product $\Lambda^k V$
admits an alternative description, as the linear space spanned by
the wedge products $v_1 \wedge \dots \wedge v_k$ of vectors $v_1,
v_2, \dots, v_k$ in $V$. Assuming $V$ comes with an inner product,
we can endow $\Lambda^k V$ with a inner product such that $\|v_1
\wedge \dots \wedge v_k\|$ is just the volume of the
$k$-dimensional parallelepiped determined by the vectors $v_1,
v_2, \dots, v_k$ in $V$.

A linear isomorphism $A: V \to W$ induces another, $\Lambda^k
A:\Lambda^k V\to \Lambda^k W$, through
$$
\Lambda^k A (v_1 \wedge \dots \wedge v_k) = Av_1 \wedge \dots
\wedge Av_k.
$$
When $V=W$, the eigenvalues of $\Lambda^k A$ are just the products
of $k$ distinct eigenvalues of $A$ (where an eigenvalue with
multiplicity $m$ is counted as $m$ ``distinct'' eigenvalues).
Correspondingly, there is a simple relation between the Lyapunov
spectra of $\Lambda^k Df$ and $Df$: the Lyapunov exponents of
$\Lambda^k Df$ are the sums of $k$ distinct Lyapunov exponents of
$Df$, with the same convention as before concerning multiplicities.
Thus,
$$
\lambda_{i_1}(x) + \lambda_{i_2}(x) + \dots + \lambda_{i_k}(x) \le
\log C_k(f)
$$
for any $1 \le i_1 < i_2 < \dots < i_k \leq d$, and our
hypothesis~\eqref{eq.condicao} implies that these sums are
strictly smaller than $\log p$, for all $k<d$.

\begin{lemma}\label{l.entropiagrande}
If $\mu$ is an invariant probability with some integrated Lyapunov
exponent less than
\begin{equation}\label{eq.definicaodec}
c(f)= \log p - \max_{1\leq k < d} \log C_k(f),
\end{equation}
then $h_\mu(f)<\log p$.
\end{lemma}

\begin{proof}
Let $\mu$ be an invariant probability, and suppose $\int
\lambda_d(x) d\mu < c(f)$. As we have just seen,
\eqref{eq.condicao} implies that $\sum_{1\leq i\leq k}
\lambda_i(x) \le \log C_k(f)$ for all $1\le k <d$. Then, using the
Ruelle inequality~\cite{Rue78},
$$
h_\mu(f) \leq \int \sum_{i: \lambda_i(x)>0} \lambda_i(x) \, d\mu
 < c(f) + \max_{1 \leq k < d} \log C_k(f) \le \log p.
$$
This proves the lemma.
\end{proof}

\section{Hyperbolic times}\label{s.hyperbolic}

For the next step we need the notion of hyperbolic times,
introduced by Alves {\it et al}~\cite{Al00,ABV00}. Given $c>0$, we
say that $n\in \natural$ is a \emph{$c$-hyperbolic time}  for
$x\in M$ if
$$
\prod _{k=0}^{j-1}  \|Df(f^{n-k}(x))^{-1}\| \leq e^{-2cj}
\quad\text{for every  $1\leq j\leq n$.}
$$
In what follows we fix $c=c(f)/10$ and speak, simply, of
\emph{hyperbolic times}. We say that $f$ has \emph{positive
density of hyperbolic times} for $x$ if the set $H_x$ of
integers which are hyperbolic times of $f$ for $x$ satisfies
\begin{equation}\label{eq.density}
\liminf_n \frac{1}{n} \#(H_x \cap[1, n]) > 0.
\end{equation}
We quote a few basic properties from \cite{ABV00}
(alternatively, see \cite{Ol03}):

%

\begin{lemma}\label{l.integral}
If a point $x$ satisfies
$$
\limsup_{n\to \infty} \frac{1}{n} \sum_{i=0}^{n-1}
\log\|Df(f^{i}(y))^{-1}\|<-4c<0,
$$
then $f$ has positive density of hyperbolic times for $x$.
\end{lemma}

In fact, the density, that is, the $\liminf$ in \eqref{eq.density},
is bounded below by some positive constant that depends only on $f$
(and our choice of $c$).

\begin{lemma}\label{l.contracao}
There exists $\delta_0>0$, depending only on $f$ and $c$, such
that given any hyperbolic time $n\ge 1$ for a point $x\in M$, and
given any $1\le j\le n$, the inverse branch $f^{-j}_{x,n}$ of
$f^j$ that sends $f^n(x)$ to $f^{n-j}(x)$ is defined on the whole
ball of radius $\delta_0$ around $f^n(x)$, and satisfies
$$
d(f^{-j}_{x,n}(z),f^{-j}_{x,n}(w))\leq e^{-jc} d(z,w)
$$
for every $z$, $w$ in that ball.
\end{lemma}

In view of Lemma~\ref{l.entropiagrande}, the next lemma applies to
any invariant measure $\mu$ with $h_\mu(f)\ge \log p$.

\begin{lemma}\label{l.iterados}
Given an invariant ergodic measure $\mu$ whose Lyapunov exponents
are all bigger than $8c$, there exists $N \in \natural$ such that
$f^N$ has positive density of hyperbolic times for $\mu$-almost
every point.
\end{lemma}

\begin{proof}
Since all Lyapunov exponents of $\mu$ are greater than $8c$, for
almost every $x\in M$ there exists $n_0(x)\ge 1$ such that
$$
\|Df^n(x)w\|\geq e^{6cn}\|w\|, \text{ for all } w\in T_xM \text{
and } n\geq n_0(x).
$$
In other words,
$$
\|Df^n(x)^{-1}\|\leq e^{-6cn},\text{ for every } n\ge n_0(x).
$$
Define $\alpha_n = \mu(\{x:n_0(x)>n\})$. Since $f$ is a local
diffeomorphism, we may also fix a constant $K>0$ such
$\|Df(x)^{-1}\|\le K$ for all $x\in M$. Then
$$
\int_M \log \|Df^n(x)^{-1}\| d\mu \leq - 6 c n + K n \alpha_n =
-(6c+K\alpha_n) n.
$$
Since $\alpha_n$ goes to zero when $n$ goes to infinity, by
choosing $N$ big enough we ensure that
$$
\int_M \frac{1}{N} \log \|Df^N(x)^{-1}\| d\mu < - 4c < 0.
$$
Then, since $\mu$ is ergodic,
$$
\lim_{n\to\infty} \frac{1}{n}\sum_{i=0}^{n-1}\frac{1}{N} \log
\|Df^N(f^{Ni}(y))^{-1}\| = \int_M \frac{1}{N} \log
\|Df^N(x)^{-1}\| d\mu < -4c.
$$
This means that we may apply Lemma~\ref{l.integral} to conclude.
\end{proof}

According to the remark following Lemma~\ref{l.integral}, we
even have that the density of hyperbolic times is bounded below
by some positive constant that depends only on $f^N$ (and our choice
of $c$).

\begin{lemma}\label{l.fubini}
Let $B\subset M$, $\theta>0$, and $g:M\to M$ be a local diffeomorphism
such that $g$ has density $>2\theta$ of hyperbolic times for every $x\in B$.
Then, given any probability measure $\nu$ on $B$ and any $m\ge 1$,
there exists $n > m$ such that
$$
\nu\big(\{x\in B: \text{$n$ is a hyperbolic time of $g$ for $x$}\}\big)
> \theta.
$$
\end{lemma}

\begin{proof}
Define $H$ to be the set of pairs $(x,n)\in B\times\natural$ such that
$n$ is a hyperbolic time for $x$.
For each $k \ge 1 $, let $\chi_k$ be the normalized counting measure
on the time interval $[m+1,m+k]$. The hypothesis implies that, given any
$x\in B$, we have
$$
\chi_k\big(H\cap(\{x\}\times\natural)\big) > 2 \theta
$$
for every sufficiently large $k$. Fix $k\ge 1$ large enough so that
this holds for a subset $C$ of points $x\in B$ with $\nu(C)>1/2$. Then,
by Fubini's theorem, $(\nu\times\chi_k)(H) > 2\theta$, and this implies
that
$$
\nu\big(H\cap(B\times\{n\})) > \theta
$$
for some $n\in[m+1,m+k]$. This gives the conclusion of the lemma.
\end{proof}

\section{Generating partitions}\label{s.generating}

In all that follows the constant $\delta_0>0$ is fixed as given by
Lemma~\ref{l.contracao}. Given a partition $\alpha$ of $M$,
we define
$$
\alpha_n=\bigvee_{j=0}^{n-1} f^{-j}(\alpha)
\quad\text{for each $n\ge 1$.}
$$

\begin{lemma}\label{l.partition}
If $\mu$ is an invariant measure such that all its Lyapunov
exponents are bigger than $8c$, and $\alpha$ is a partition with
diameter less than $\delta_0$, for $\mu$-almost every $x\in M$,
the diameter of $\alpha_n(x)$ goes to zero when $n$ goes to $\infty$.
In particular, $\alpha$ is an $f$-generating partition with respect
to $\mu$.
\end{lemma}

\begin{proof}
By Lemma~\ref{l.iterados} there exists $N\ge 1$ such that $f^N$
has positive density of hyperbolic times for $\mu$-almost every
point. Define
$$
\gamma_k=\bigvee_{j=0}^{k-1} f^{-jN}(\alpha)
 \quad\text{for each $k\ge 1$.}
$$
By Lemma~\ref{l.contracao}, if $k$ is a hyperbolic time of $f^N$ for
$x$ then $\diam\gamma_k(x)\le e^{-cn}$. In particular, since the sets
$\gamma_k(x)$ are non-increasing with $k$, the diameter of $\gamma_k(x)$
goes to zero when $k\to\infty$.
Since $\alpha_{kN}(x)\subset \gamma_k(x)$ and the sequence
$\diam\alpha_n(x)$ is non-increasing, this immediatelly gives that
the diameter of $\alpha_n(x)$ goes to zero when $n$ goes to infinity,
for $\mu$-almost every $x\in M$.

The rest of the argument is very standard. It goes as follows. To
prove that $\alpha$ is a generating partition for $f$ with respect
to $\mu$, it suffices to show that, given any measurable set $E$ and
any $\vep>0$, there exists $n\ge 1$ and elements
$A_n^i$, $i=1, \dots, m(n)$ of $\alpha_n$ such that
$$
\mu\big(\bigcup_{i=1}^m A^i_n \Delta E\big) < \vep.
$$
Consider compact sets $K_1\subset A$ and $K_2\subset A^c$ such that
$\mu(K_1 \Delta A)$ and $\mu(K_2 \Delta A^c)$ are both less than
$\vep/4$. Fix $n\ge 1$ large enough so that $\diam\alpha_n(x)$ is
smaller than the distance from $K_1$ to $K_2$ outside a set of points
$x$ with measure less than $\vep/4$. Let $A_n^i$, $i=1, \ldots, m(n)$
be the sets $\alpha_n(x)$ that intersect $K_1$. Then, they are
all disjoint from $K_2$, and so $\mu(\bigcup_i A^i_n \Delta E)$
is bounded above by
$$
\mu(E\setminus \bigcup_i A_n^i)+\mu(\bigcup_i A_n^i \setminus E)
\leq \mu(E \setminus K_1) + \mu(E^c \setminus K_2) + \vep/4 \leq \vep.
$$
This completes the proof.
\end{proof}

\section{Rokhlin's formula}\label{s.rokhlin}

The \emph{Jacobian} of a measure $\mu$ with respect to $f$ is the
(essentially unique) function $J_\mu f$ satisfying
$$
\mu(f(A))=\int_A J_\mu f d\mu.
$$
for any measurable set $A$ such that $f|_A$ is injective.

In other words, the Jacobian is defined by $J_\mu f = {d(\mu \circ f)/d\mu}$. Jacobians for every measure do
exist in this context, because $f$ is finite-to-one (countable-to-one would suffice). Using the definition, one
can verify that $J_\mu f^n(x) = \prod_{i=0}^{n-1} J_\mu f(f^i(x))$ is a Jacobian for each $f^n$. In the case of
$\mu$ is an invariant measure, we observe that from the definition follows that $J_\mu f \geq 1$ in $\mu$-almost
everywhere.

Let $f:M\to M$ be a measurable transformation, $\mu$ be an
invariant probability. Suppose there exists a finite or countable
partition $\alpha$ of $M$ such that
\begin{enumerate} \label{e.conditions}
\item[(a)] $f$ is locally injective, meaning that it is injective
on every atom of $\alpha$;
\item[(b)] $\alpha$ is $f$-generating with respect to $\mu$, in
the sense that $\diam\alpha_n(x)\to 0$ for $\mu$-almost every $x$.
\end{enumerate}

\begin{proposition}\label{p.rokhlinformula}
If $\mu$ is an invariant measure satisfying (a) and (b) as above,
then
$$
h_\mu(f)=\int\log J_\mu f\,d\mu,
$$ where $J_\mu f$ denotes any
Jacobian of $f$ relative to $\mu$.
\end{proposition}

Let $\alpha_\infty=\bigvee_{j=0}^\infty f^{-j}(\alpha)$. Denote
$\beta_\infty=\bigvee_{j=1}^\infty f^{-j}(\alpha)$ and
$\beta_n=\bigvee_{j=1}^n f^{-j}(\alpha)$ for each $n\ge 1$. Notice
that $\beta_\infty(x)=f^{-1}(\alpha_\infty(f(x)))$. The hypothesis
that $\alpha$ is generating implies that $\alpha_\infty(x)=\{x\}$,
and so
\begin{equation}\label{eq.sat}
\beta_\infty(x)=\{f^{-1}(f(x))\} \quad\text{for $\mu$-almost all $x\in
M$}.
\end{equation}

The conditional expectation of a function $\varphi:M\to\real$
relative to a partition $\gamma$ is the essentially unique
$\gamma$-measurable function $E_\mu(\varphi \mid \gamma)$ such that
\begin{equation}\label{eq.expect1}
\int_B E_\mu(\varphi \mid \gamma) \, d\mu = \int_B \varphi\,d\mu
\end{equation}
for every $\gamma$-measurable set $B$.

\begin{lemma}\label{l.Emu}
$E_\mu(\varphi \mid \beta_\infty)(x) = \sum_{y\in\beta_\infty(x)}
\frac{1}{J_\mu f(y)} \varphi(y)$ for $\mu$-almost every $x$.
\end{lemma}

\begin{proof}
It is clear that the function on the right hand side is
$\beta_\infty$-measurable. Let $B$ be any
$\beta_\infty$-measurable set, that is, any measurable set that
consists of entire atoms of $\beta_\infty$. By \eqref{eq.sat},
there exists a measurable set $C$ such that $B=f^{-1}(C)$. Then,
since $\mu$ is invariant,
$$
\begin{aligned}
\int_B \sum_{y\in\beta_\infty(x)}\frac{1}{J_\mu f(y)}\varphi(y)\,
d\mu(x)
 & = \int_C \sum_{y\in f^{-1}(z)}\frac{1}{J_\mu f(y)}\varphi(y)\, d\mu(z)
 \\
 & = \sum_{A\in\alpha} \int_{C_A} \frac{1}{J_\mu f(y_A)}\varphi(y_A)\, d\mu(z),
\end{aligned}
$$ where $C_A=f(B\cap A)$ and $y_A=(f\mid A)^{-1}(z)$. Since
every $f\mid A$ is injective, we may use the definition of the
Jacobian to rewrite the latter expression as
$$
\sum_{A\in\alpha} \int_{B\cap A} \varphi(y)\, d\mu(y)
 = \int_B \varphi\, d\mu.
$$
This proves \eqref{eq.expect1} and the lemma.
\end{proof}

For $*\in\natural \cup\{\infty\}$, define the conditional entropy
(Definition~4.8 in~\cite{Wa82})
\begin{equation}\label{eq.eq.relativa}
H_\mu(\alpha \mid\beta_*) = \int \sum_{A\in\alpha}
-E_\mu(\chi_A\mid\beta_*) \log E_\mu(\chi_A\mid\beta_*) \,d\mu.
\end{equation}

\begin{lemma}\label{l.entropia}~
\begin{enumerate}
\item $H_\mu(\alpha \mid \beta_n)
 = \sum_{A\in\alpha} \sum_{B\in\beta_n} - \mu(A\cap B)\log\frac{\mu(A\cap B)}{\mu(B)}$
 for $n\in\natural$.
\item $H_\mu(\alpha \mid \beta_\infty) = \int \log J_\mu f\,d\mu$.
\end{enumerate}
\end{lemma}

\begin{proof}
For $n\in\natural$, the partition $\beta_n$ is countable, and so
$$
E_\mu(\chi_A \mid \beta_n)(x)
 = \frac{1}{\mu(\beta_n(x))} \int_{\beta_n(x)} \chi_A \,d\mu
 = \frac{\mu(A\cap\beta_n(x))}{\mu(\beta_n(x))}
$$
for every $A\in\alpha$. It follows that
$$
H_\mu(\alpha \mid \beta_n)
 = \sum_{A\in\alpha} \sum_{B\in\beta_n} \int_B
 - \frac{\mu(A\cap\beta_n(x))}{\mu(\beta_n(x))} \log\frac{\mu(A\cap\beta_n(x))}{\mu(\beta_n(x))}
 \,d\mu(x).
$$
This gives the first statement. Next, Lemma~\ref{l.Emu} says that
$$
E_\mu(\chi_A \mid \beta_\infty)=\psi_A\circ f, \quad\text{where }
\psi_A(z) = \sum_{y\in f^{-1}(z)} \frac{1}{J_\mu f(y)} \chi_A(y).
$$
Notice that if $z\in f(A)$ then $\psi_A(z)=1/J_\mu f(y_A)$, where
$y_A=(f\mid A)^{-1}(z)$, and if $z\notin f(A)$ then $\psi_A(z)=0$.
Therefore,
$$
\begin{aligned}
H_\mu(\alpha\mid\beta_\infty)
 & = \int \sum_{A\in\alpha} -\psi_A(z)\log\psi_A(z)\,d\mu(z)
 \\ & = \sum_{A\in\alpha} \int_{f(A)} \frac{1}{J_\mu f(y_A)} \log J_\mu
 f(y_A)\,d\mu(z).
\end{aligned}
$$
Using the definition of Jacobian, and the assumption that $f$ is
injective on $A$, this gives
$$
H_\mu(\alpha\mid\beta_\infty)
 = \sum_{A\in\alpha} \int_{A}\log J_\mu f(y)\,d\mu(y)
 = \int \log J_\mu f(y)\,d\mu(y),
$$
as claimed.
\end{proof}

\begin{proof}[Proof of Proposition~\ref{p.rokhlinformula}]
Since the partition $\alpha$ is generating, $h_\mu(f) =
h_\mu(f,\alpha)$. Then,
$$
h_\mu(f,\alpha)
 = \lim_n H_\mu(\alpha \mid\beta_n)
 = H_\mu(\alpha \mid\beta_\infty),
$$
by Theorem~4.14 of~\cite{Wa82}. Combined with the second part of
Lemma~\ref{l.entropia}, this gives $h_\mu(f) = \int \log J_\mu f
\, d\mu$, as claimed.
\end{proof}

\section{Existence}\label{s.existence}

Here we prove that every maximal eigenmeasure is a maximizing
measure. The first step is

\begin{lemma}\label{l.jacobian}
If $\mu$ is a maximal eigenmeasure then $J_\mu f$ is constant
equal to $p$.
\end{lemma}

\begin{proof}
Let $A$ be any measurable set such that $f|_A$ is injective. Take
a sequence $\{g_n\} \in C(M)$ such that $g_n\rightarrow \chi_A$ at
$\mu$-almost every point and $\sup |g_n| \leq 2$ for all $n$. By
definition,
$$
\SL g_n (x) =\sum_{f(y)=x} g_n(y).
$$
The last expression converges to $\chi_{f(A)}(x)$ at $\mu$-almost
every point. Hence, by the dominated convergence theorem,
$$
\int p \, g_n \, d\mu
 = \int g_n \, d(\SL^*\mu)
 = \int \SL g_n \, d\mu
 \to \mu(f(A)).
$$
Since the left hand side also converges to $\int_A p d\mu$, we
conclude that
$$
\mu(f(A))=\int_A p d\mu,
$$
which proves the lemma.
\end{proof}

\begin{lemma}\label{l.eigenentropy}
If $\mu$ is a maximal eigenmeasure then $h_\mu(f)\ge\log p$.
\end{lemma}

\begin{proof}
We define the \emph{dynamical ball}  $\SB_\epsilon(n,x)$ by
$$
\SB_\epsilon(n,x)=\{y \in M; d(f^i(x),f^i(y))< \epsilon, \text{ for } i=0,\dots,n-1\}.
$$
If $\epsilon$ small enough so that $f^n|_{\SB_\epsilon(n,x)}$ is
injective, then:
$$
1 = \mu(M) \geq \mu(f^n(\SB_\epsilon(n,x))) = p^n
\mu(\SB_\epsilon(n,x)).
$$
In particular, we may conclude that
$$
- \limsup (1/n) \log \mu(\SB_\epsilon(n,x)) \geq \log p
$$
for every $n$ and $\epsilon$ small. By the Brin-Katok local
entropy formula (see \cite{Man87})
$$
h_\mu(f) = -\int \lim\limits_{\epsilon \rightarrow 0}
                  \lim\limits_{n\rightarrow \infty}
 \frac{1}{n} \log \mu(\SB_\epsilon(n,x)) \, d\mu(x) \geq \log p.
$$
This proves the lemma.
\end{proof}

\begin{corollary}
Every maximal eigenmeasure $\mu$ has entropy equal to $\log p$.
\end{corollary}

\begin{proof}
By Lemma~\ref{l.eigenentropy}, the entropy is at least $\log p$.
Then, we may apply Lemma~\ref{l.entropiagrande} to conclude that
all Lyapunov exponents of $\mu$ are positive. It follows, by
Lemma~\ref{l.partition}, that $\mu$ admits generating partitions
with small diameter. Hence, we may apply
Proposition~\ref{p.rokhlinformula} and Lemma~\ref{l.jacobian}, to
find that $h_{\mu}(f) = \int \log J_\mu d\mu = \log p$.
\end{proof}

\begin{lemma}\label{l.grau}
The topological entropy $h_{top}(f) = \log p$. Moreover, if $\eta$
is any ergodic maximizing measure then the Jacobian $J_\eta f$ is
constant equal to $p$.
\end{lemma}

\begin{proof}
Consider any probability $\eta$ such $h_\eta(f) \geq \log p$. By
Lemma~\ref{l.entropiagrande} all Lyapunov exponents of $\eta$ are
bigger than $c(f)$. Then, by Lemma~\ref{l.partition}, there exist
generating partitions with arbitrarily small diameter. This
ensures we may apply Proposition~\ref{p.rokhlinformula} to $\eta$.
We get that
$$
h_\eta(f) = \int \log J_\eta f \, d\eta.
$$
Let us write $g_\eta = 1/(J_\eta f)$. The assumption that $\eta$
is invariant means that
$$
\sum_{f(y)=x} g_\eta(y)=1
$$
for $\eta$-almost every $x\in M$. From the previous equality, we
find
\begin{equation}\label{eq.quaseentropia}
0 \leq h_\eta(f)-\log p
  = \int \log \frac{p^{-1}}{g_\eta}\,d\eta
  = \int \sum_{f(y)=x} g_\eta(y)\log\frac{p^{-1}}{g_\eta(y)}\,d\eta,
\end{equation}
where the first equality uses $g_\eta=1/J_\eta f$. Using the Jensen inequality:


$$
\sum_{f(y)=x}  g_\eta(y)\log \frac{p^{-1}}{g_\eta(y)}
 \le \log\big(\sum_{f(y)=x} g_\eta(y) \frac{p^{-1}}{g_\eta(y)}\big)
 = \log\big(\sum_{f(y)=x} p^{-1}\big)=0
$$
at $\eta$-almost every point. Since the integral is non-negative,
by \eqref{eq.quaseentropia}, the equality must hold $\eta$-almost
everywhere, and $h_\eta(f) - \log p = 0$. Since $\eta$ is
arbitrary, this proves that $\log p = h_{top}(f)$.

\relax From the last part of Lemma~\ref{l.calculus}, we get that the values of $\log{p^{-1}}/{g_\eta(y)}$ are
the same for all $y\in f^{-1}(x)$. In other words, for $\eta$-almost every $x\in M$ there exists a number $c(x)$
such that ${p^{-1}}/{g_\eta(y)} = c(x)$ for every $y\in f^{-1}(x)$. Then
$$
\frac{1}{c(x)}
 = \sum_{y\in f^{-1}(x)} \frac{p^{-1}}{c(x)}
 = \sum_{y\in f^{-1}(x)} g_\eta(y)
 = 1
$$
for $\eta$-almost every $x$. This means, precisely, that $J_\eta
f(y) = p$ for every $y$ on the pre-image of a full $\eta$-measure
set.
\end{proof}

\section{Uniqueness}\label{s.uniqueness}

In this section we assume $f$ is topologically mixing, and conclude
that the maximizing measure is unique and supported on the whole
ambient $M$. It suffices to consider ergodic measures, because the
ergodic components of maximizing measures are also maximizing measures.

\begin{lemma}\label{l.support}
Any ergodic maximizing measure $\mu$ is supported on the whole $M$.
\end{lemma}

\begin{proof}
Suppose $\mu(U)=0$ for some non-empty open set $U$. By the mixing assumption, there exists $N\ge 1$ such
$f^N(U)=M$. Partitioning $U$ into subsets $U_1$, \ldots, $U_k$ such that every $f^N|U_j$ is injective, we get
that
$$
\mu(f^N(U_j))=\int_{U_j} J_\mu f^N \,d\mu = 0
$$
for $j=1, \dots, k$. Recall Lemma~\ref{l.grau}. This implies that $\mu(M)=\mu((f^N(U))=0$, which is a
contradiction.
\end{proof}

This has the following useful consequence: given any $\delta>0$ there
exist $b=b(\delta)>0$ such that
\begin{equation}\label{eq.bedelta}
\mu(B(x,\delta)) \ge b
\quad\text{for all } x\in M.
\end{equation}
Indeed, if there were points such that the balls of radius $\delta$
around them have arbitrarily small measures then, considering an
accumulation point, one would get a ball with zero measure, and that
would contradict Lemma~\ref{l.support}.

Now let $\mu_1$ and $\mu_2$ be any two ergodic maximizing measures.
Our goal is to prove that the two measures coincide. As a first
step we prove that they are equivalent. For this, we fix any
(finite) partition $\SP$ of $M$ into subsets $P$ such that $P$ has
non-empty interior, and the boundary $\partial P$ has zero measure
for both $\mu_1$ and $\mu_2$. Fixing $\delta>0$ small so that every
$P\in\SP$ contains some ball of radius $\delta$, and applying
\eqref{eq.bedelta} to both measures, we conclude that there exists
$B>0$ such that
\begin{equation}\label{eq.equivalencia1}
\mu_1(P)\le B\mu_2(P)
\text{\ and\ }
\mu_2(P)\le B\mu_1(P)
\quad\text{for all $P\in\SP$.}
\end{equation}
Now let $g$ be an inverse branch of any iterate $f^n$, $n\ge 1$.
Using Lemma~\ref{l.grau}, we get that
$\mu_i(P)=p^n\mu_i(g(P))$ for $i=1, 2$. It follows that
\eqref{eq.equivalencia1} remains valid for the images $g(P)$:
\begin{equation}\label{eq.equivalencia2}
\mu_1(g(P))\le B\mu_2(g(P))
\text{\ and\ }
\mu_2(g(P)) \le B \mu_1(g(P))
\end{equation}
for every $P\in\SP$ and every inverse branch $g$ of $f^n$,
for any $n\ge 1$. We denote by $\SQ$ the family of all
such images $g(P)$.

\begin{lemma}\label{l.exaustao}
Given any measurable set $E\subset M$ and any $\vep>0$ there exists
a family $\SE$ of pairwise disjoint elements of $\SQ$ such that
$$
\mu_i\big(E \setminus \bigcup_{\SE} g(P)\big) = 0
\quad\text{and}\quad
\mu_i\big(\bigcup_{\SE} g(P) \setminus E\big)\le\vep
\text{ for $i=1, 2$.}
$$
\end{lemma}

\begin{proof}
By Lemma~\ref{l.entropiagrande}, all Lyapunov exponents of $\mu_i$
are larger than $c(f)$. Hence, by Lemma~\ref{l.iterados} and the
remark following it, there exists $N\ge 1$ and $\theta>0$ such that
$\mu_i$-almost every point has density $> 2\theta$ of hyperbolic
times.

Let $U_1$ be an open set and $K_1$ be a compact set such that
$K_1\subset E\subset U_1$ and $\mu_i(U_1\setminus E)\le\vep$
for $i=1, 2$ and $\mu_i(K_1)\ge(1/2)\mu(E)$.
Using Lemma~\ref{l.fubini} with $B=K_1$ and $\nu=\mu_i/\mu_i(K_1)$,
we may find $n_1\ge 1$ such that $e^{-c n_1} < d(K_1,U_1^c)$
and the subset $L_1$ of points $x\in K_1$ for which $n_1$ is a
hyperbolic time satisfies
$\mu_i(L_1) \ge \theta \mu_i(K_1) \ge (\theta/2)\mu_i(E)$.
Let $\SE_1$ the family of all $g(P)$ that intersect $L_1$, with
$P\in\SP$ and $g$ an inverse branch of $f^{n_1}$.
Notice that the elements of $\SE_1$ are pairwise disjoint,
because the elements of $P$ are pairwise disjoint.
Moreover, by Lemma~\ref{l.contracao}, their diameter is less
than $e^{-c n_1}$. Thus, the union $E_1$ of all the elements
of $\SE_1$ is contained in $U_1$. By construction, it satisfies
$$
\mu_i(E_1\cap E) \geq \mu_i(L_1)
\ge \theta \mu_i(K_1) \ge (\theta/2) \mu_i(E).
$$

Next, consider the open set $U_2=U_1\setminus \overline E_1$
and let $K_2 \subset E\setminus\overline E_1$ be a compact set
such that $\mu_i(K_2)\ge (1/2) \mu_i(E\setminus E_1)$.
Observe $\mu_i(\overline E_1\setminus E_1)=0$ because the
boundaries of the atoms of $\SP$ have zero measure and that
is preserved by the inverse branches, since $\mu_i$ is invariant.
Reasoning as before, we may find $n_2>n_1$ such that
$e^{-c n_2} < d(K_2,U_2^c)$ and a set $L_2\subset K_2$ such that
$\mu_i(L_2)\ge\theta\mu_i(K_2)$ and $n_2$ is a hyperbolic time
for every $x\in L_2$. Denote by $\SE_2$ the family of inverse
images $g(P)$ that intersect $L_2$, with $P\in\SP$ and $g$ an
inverse branch of $f^{n_2}$. As before, the elements of $\SE_2$
are pairwise disjoint, and their diameters are smaller than
$e^{-c n_2}$. The latter ensures that their union $E_2$ is
contained in $U_2$. Consequently, the elements of the union
$\SE_1\cup\SE_2$ are also pairwise disjoint. Moreover,
$$
\mu_i\big(E_2\cap [E\setminus E_1]\big)
\ge \mu_i(L_2) \ge \theta \mu_i(K_2)
\ge (\theta/2) \mu_i(E \setminus E_1).
$$

Repeating this procedure, we construct families $\SE_k$, $k\ge 1$
of elements of $\SQ$ such that their elements are all pairwise
disjoint and contained in $U_1$, and
\begin{equation}\label{eq.interests}
\mu_i\big(E_{k+1}\cap [E\setminus (E_1\cup\cdots\cup E_k)]\big)
\ge (\theta/2) \mu_i\big(E\setminus (E_1\cup\cdots\cup E_k)\big)
\end{equation}
for all $k\ge 1$, where $E_j=\cup_{\SE_j} g(P)$. Thus,
$\mu_i\big(\bigcup_{k=1}^\infty E_k \setminus E\big)
\le \mu_i(U_1\setminus E) \le \vep$ for $i=1, 2$, and
\eqref{eq.interests} implies that
$$
\mu_i(E \setminus \bigcup_{k=1}^\infty E_k)=0.
$$
This completes the proof of the lemma,
with $\SE=\bigcup_{k=1}^\infty \SE_k$.
\end{proof}

\begin{remark}
The lemma remains true if one asks that
$ \mu_i\big(E \setminus \bigcup_{\SE} g(P)\big) = 0$ for both $i=1, 2$.
This follows from a variation of the previous construction, considering
each one of the two measures alternately: for each $k\ge 1$ consider
$i\equiv k \mod 2$; then ask that
$\mu_i(K_k)\ge (1/2) \mu_i(E\setminus[E_1\cup\cdots\cup E_k])$,
and choose $n_k$ such that $\mu_i(L_k)\ge\theta\mu(K_k)$.
The same kind of argument applies with any number of probability measures
$\mu_1, \ldots, \mu_r$. These extensions will not be used here.
\end{remark}

Combining \eqref{eq.equivalencia2} with Lemma~\ref{l.exaustao},
we get that, for any measurable set $E\subset M$,
$$
\mu_1(E)\le \vep + \sum_{\SE} \mu_1(g(P))
\le \vep + B \sum_{\SE} \mu_2(g(P)) = \vep + B \mu_2(E).
$$
As $\vep>0$ is arbitrary, we get that $\mu_1(E)\le B\mu_2(E)$.
A symmetric argument gives that $\mu_2(E)\le B\mu_1(E)$ for any
measurable set $E$.
This implies that $\mu_1 = h \mu_2$ where the Radon-Nikodym derivative
$h$ satisfies $B^{-1}\leq h \leq B$. Since $\mu_1$ and $\mu_2$ are
invariant measures,
$$
\mu_1 = f_* \mu_1 = (h\circ f) f_* \mu_2 = (h \circ f) \mu_2.
$$
As the Radon-Nikodym derivative is essentially unique, we get that
$h = h\circ f$ at $\mu_2$-almost every point. By ergodicity, it
follows that $h$ is constant almost everywhere. Since the $\mu_i$ are
both probabilities, we get that $h=1$ and so $\mu_1=\mu_2$.
This proves uniqueness of the maximizing measure.

\vspace{1cm}

\noindent Krerley Oliveira ( krerley{\@@}mat.ufal.br )\\
Departamento de Matem\'atica - UFAL, Campus A.C. Sim\~oes, s/n
57072-090 Macei\'o, Alagoas - Brazil

\bigskip

\noindent Marcelo Viana ( viana\@@impa.br ) \\
IMPA, Est. D. Castorina 110 \\
22460-320 Rio de Janeiro, RJ, Brazil

\end{document}